\documentclass[10.5pt]{article}
\usepackage{amsmath}
\usepackage{amssymb}
\usepackage{amscd}
\usepackage{xspace}
\usepackage{verbatim}
\newcommand{\mysection}[1]{
\section{#1}\setcounter{equation}{0}}
\title{\bf Boundary singularities of $N$-harmonic functions
\footnote {To appear in{ \it Communications in Partial Differential Equations}}}
\author{{\bf Rouba Borghol}, {\bf Laurent V\'eron}\\
{\small Department of Mathematics,}\\
 {\small  University of Tours,  FRANCE}
}

\date{}
\begin{document}
\maketitle

\newcommand{\txt}[1]{\;\text{ #1 }\;}
\newcommand{\tbf}{\textbf}
\newcommand{\tit}{\textit}
\newcommand{\tsc}{\textsc}
\newcommand{\trm}{\textrm}
\newcommand{\mbf}{\mathbf}
\newcommand{\mrm}{\mathrm}
\newcommand{\bsym}{\boldsymbol}
\newcommand{\scs}{\scriptstyle}
\newcommand{\sss}{\scriptscriptstyle}
\newcommand{\txts}{\textstyle}
\newcommand{\dsps}{\displaystyle}
\newcommand{\fnz}{\footnotesize}
\newcommand{\scz}{\scriptsize}
\newcommand{\be}{
\begin{equation}
}
\newcommand{\bel}[1]{
\begin{equation}
\label{#1}}
\newcommand{\ee}{
\end{equation}
}
\newcommand{\eqnl}[2]{
\begin{equation}
\label{#1}{#2}
\end{equation}
}
\newtheorem{subn}{\name}
\renewcommand{\thesubn}{}
\newcommand{\bsn}[1]{\def\name{#1}
\begin{subn}}
\newcommand{\esn}{
\end{subn}}
\newtheorem{sub}{\name}[section]
\newcommand{\dn}[1]{\def\name{#1}}   
\newcommand{\bs}{
\begin{sub}}
\newcommand{\es}{
\end{sub}}
\newcommand{\bsl}[1]{
\begin{sub}\label{#1}}
\newcommand{\bth}[1]{\def\name{Theorem}
\begin{sub}\label{t:#1}}
\newcommand{\blemma}[1]{\def\name{Lemma}
\begin{sub}\label{l:#1}}
\newcommand{\bcor}[1]{\def\name{Corollary}
\begin{sub}\label{c:#1}}
\newcommand{\bdef}[1]{\def\name{Definition}
\begin{sub}\label{d:#1}}
\newcommand{\bprop}[1]{\def\name{Proposition}
\begin{sub}\label{p:#1}}
\newcommand{\R}{\eqref}
\newcommand{\rth}[1]{Theorem~\ref{t:#1}}
\newcommand{\rlemma}[1]{Lemma~\ref{l:#1}}
\newcommand{\rcor}[1]{Corollary~\ref{c:#1}}
\newcommand{\rdef}[1]{Definition~\ref{d:#1}}
\newcommand{\rprop}[1]{Proposition~\ref{p:#1}}
\newcommand{\BA}{
\begin{array}}
\newcommand{\EA}{
\end{array}}
\newcommand{\BAN}{\renewcommand{\arraystretch}{1.2}
\setlength{\arraycolsep}{2pt}
\begin{array}}
\newcommand{\BAV}[2]{\renewcommand{\arraystretch}{#1}
\setlength{\arraycolsep}{#2}
\begin{array}}
\newcommand{\BSA}{
\begin{subarray}}
\newcommand{\ESA}{
\end{subarray}}
\newcommand{\BAL}{
\begin{aligned}}
\newcommand{\EAL}{
\end{aligned}}
\newcommand{\BALG}{
\begin{alignat}}
\newcommand{\EALG}{
\end{alignat}}
\newcommand{\BALGN}{
\begin{alignat*}}
\newcommand{\EALGN}{
\end{alignat*}}
\newcommand{\note}[1]{\textit{#1.}\hspace{2mm}}
\newcommand{\Proof}{\note{Proof}}
\newcommand{\qeda}{\hspace{10mm}\hfill $\square$}
\newcommand{\qed}{\\
${}$ \hfill $\square$}
\newcommand{\Remark}{\note{Remark}}
\newcommand{\modin}{$\,$\\
[-4mm] \indent}
\newcommand{\forevery}{\quad \forall}
\newcommand{\set}[1]{\{#1\}}
\newcommand{\setdef}[2]{\{\,#1:\,#2\,\}}
\newcommand{\setm}[2]{\{\,#1\mid #2\,\}}
\newcommand{\lra}{\longrightarrow}
\newcommand{\lla}{\longleftarrow}
\newcommand{\llra}{\longleftrightarrow}
\newcommand{\Lra}{\Longrightarrow}
\newcommand{\Lla}{\Longleftarrow}
\newcommand{\Llra}{\Longleftrightarrow}
\newcommand{\warrow}{\rightharpoonup}
\newcommand{
\paran}[1]{\left (#1 \right )}
\newcommand{\sqbr}[1]{\left [#1 \right ]}
\newcommand{\curlybr}[1]{\left \{#1 \right \}}
\newcommand{\abs}[1]{\left |#1\right |}
\newcommand{\norm}[1]{\left \|#1\right \|}
\newcommand{
\paranb}[1]{\big (#1 \big )}
\newcommand{\lsqbrb}[1]{\big [#1 \big ]}
\newcommand{\lcurlybrb}[1]{\big \{#1 \big \}}
\newcommand{\absb}[1]{\big |#1\big |}
\newcommand{\normb}[1]{\big \|#1\big \|}
\newcommand{
\paranB}[1]{\Big (#1 \Big )}
\newcommand{\absB}[1]{\Big |#1\Big |}
\newcommand{\normB}[1]{\Big \|#1\Big \|}

\newcommand{\thkl}{\rule[-.5mm]{.3mm}{3mm}}
\newcommand{\thknorm}[1]{\thkl #1 \thkl\,}
\newcommand{\trinorm}[1]{|\!|\!| #1 |\!|\!|\,}
\newcommand{\bang}[1]{\langle #1 \rangle}
\def\angb<#1>{\langle #1 \rangle}
\newcommand{\vstrut}[1]{\rule{0mm}{#1}}
\newcommand{\rec}[1]{\frac{1}{#1}}
\newcommand{\opname}[1]{\mbox{\rm #1}\,}
\newcommand{\supp}{\opname{supp}}
\newcommand{\dist}{\opname{dist}}
\newcommand{\myfrac}[2]{{\displaystyle \frac{#1}{#2} }}
\newcommand{\myint}[2]{{\displaystyle \int_{#1}^{#2}}}
\newcommand{\mysum}[2]{{\displaystyle \sum_{#1}^{#2}}}
\newcommand {\dint}{{\displaystyle \int\!\!\int}}
\newcommand{\q}{\quad}
\newcommand{\qq}{\qquad}
\newcommand{\hsp}[1]{\hspace{#1mm}}
\newcommand{\vsp}[1]{\vspace{#1mm}}
\newcommand{\ity}{\infty}
\newcommand{\prt}{
\partial}
\newcommand{\sms}{\setminus}
\newcommand{\ems}{\emptyset}
\newcommand{\ti}{\times}
\newcommand{\pr}{^\prime}
\newcommand{\ppr}{^{\prime\prime}}
\newcommand{\tl}{\tilde}
\newcommand{\sbs}{\subset}
\newcommand{\sbeq}{\subseteq}
\newcommand{\nind}{\noindent}
\newcommand{\ind}{\indent}
\newcommand{\ovl}{\overline}
\newcommand{\unl}{\underline}
\newcommand{\nin}{\not\in}
\newcommand{\pfrac}[2]{\genfrac{(}{)}{}{}{#1}{#2}}

\def\ga{\alpha}     \def\gb{\beta}       \def\gg{\gamma}
\def\gc{\chi}       \def\gd{\delta}      \def\ge{\epsilon}
\def\gth{\theta}                         \def\vge{\varepsilon}
\def\gf{\phi}       \def\vgf{\varphi}    \def\gh{\eta}
\def\gi{\iota}      \def\gk{\kappa}      \def\gl{\lambda}
\def\gm{\mu}        \def\gn{\nu}         \def\gp{\pi}
\def\vgp{\varpi}    \def\gr{\rho}        \def\vgr{\varrho}
\def\gs{\sigma}     \def\vgs{\varsigma}  \def\gt{\tau}
\def\gu{\upsilon}   \def\gv{\vartheta}   \def\gw{\omega}
\def\gx{\xi}        \def\gy{\psi}        \def\gz{\zeta}
\def\Gg{\Gamma}     \def\Gd{\Delta}      \def\Gf{\Phi}
\def\Gth{\Theta}
\def\Gl{\Lambda}    \def\Gs{\Sigma}      \def\Gp{\Pi}
\def\Gw{\Omega}     \def\Gx{\Xi}         \def\Gy{\Psi}

\def\CS{{\mathcal S}}   \def\CM{{\mathcal M}}   \def\CN{{\mathcal N}}
\def\CR{{\mathcal R}}   \def\CO{{\mathcal O}}   \def\CP{{\mathcal P}}
\def\CA{{\mathcal A}}   \def\CB{{\mathcal B}}   \def\CC{{\mathcal C}}
\def\CD{{\mathcal D}}   \def\CE{{\mathcal E}}   \def\CF{{\mathcal F}}
\def\CG{{\mathcal G}}   \def\CH{{\mathcal H}}   \def\CI{{\mathcal I}}
\def\CJ{{\mathcal J}}   \def\CK{{\mathcal K}}   \def\CL{{\mathcal L}}
\def\CT{{\mathcal T}}   \def\CU{{\mathcal U}}   \def\CV{{\mathcal V}}
\def\CZ{{\mathcal Z}}   \def\CX{{\mathcal X}}   \def\CY{{\mathcal Y}}
\def\CW{{\mathcal W}} \def\CQ{{\mathcal Q}} 
\def\BBA {\mathbb A}   \def\BBb {\mathbb B}    \def\BBC {\mathbb C}
\def\BBD {\mathbb D}   \def\BBE {\mathbb E}    \def\BBF {\mathbb F}
\def\BBG {\mathbb G}   \def\BBH {\mathbb H}    \def\BBI {\mathbb I}
\def\BBJ {\mathbb J}   \def\BBK {\mathbb K}    \def\BBL {\mathbb L}
\def\BBM {\mathbb M}   \def\BBN {\mathbb N}    \def\BBO {\mathbb O}
\def\BBP {\mathbb P}   \def\BBR {\mathbb R}    \def\BBS {\mathbb S}
\def\BBT {\mathbb T}   \def\BBU {\mathbb U}    \def\BBV {\mathbb V}
\def\BBW {\mathbb W}   \def\BBX {\mathbb X}    \def\BBY {\mathbb Y}
\def\BBZ {\mathbb Z}

\def\GTA {\mathfrak A}   \def\GTB {\mathfrak B}    \def\GTC {\mathfrak C}
\def\GTD {\mathfrak D}   \def\GTE {\mathfrak E}    \def\GTF {\mathfrak F}
\def\GTG {\mathfrak G}   \def\GTH {\mathfrak H}    \def\GTI {\mathfrak I}
\def\GTJ {\mathfrak J}   \def\GTK {\mathfrak K}    \def\GTL {\mathfrak L}
\def\GTM {\mathfrak M}   \def\GTN {\mathfrak N}    \def\GTO {\mathfrak O}
\def\GTP {\mathfrak P}   \def\GTR {\mathfrak R}    \def\GTS {\mathfrak S}
\def\GTT {\mathfrak T}   \def\GTU {\mathfrak U}    \def\GTV {\mathfrak V}
\def\GTW {\mathfrak W}   \def\GTX {\mathfrak X}    \def\GTY {\mathfrak Y}
\def\GTZ {\mathfrak Z}   \def\GTQ {\mathfrak Q}

\font\Sym= msam10 
\def\SYM#1{\hbox{\Sym #1}}
\newcommand{\bdw}{\prt\Gw\xspace}
\medskip
\mysection {Introduction}
Let $\Gw$ be a domain is $\BBR^N$ ($N\geq 2$) with a $C^2$ compact boundary $\prt\Gw$. A function $u\in W_{loc}^{1,p}(\Gw)$ is $p$-harmonic if
\begin{equation}\label {p-harm}
\int_\Gw\abs {Du}^{p-2}\langle Du,D\gf\rangle\,dx=0
\end {equation}
for any $\gf\in C^1_0(\Gw)$. Such functions are locally $C^{1,\ga}$ for some $\ga\in (0,1)$. In the case $p=N$, the function $u$ is called $N$-harmonic. The $N$-harmonic functions play an important role as a natural extension of classical harmonic functions. They also appear in the theory of bounded distortion mappings \cite {Re}. One of the main properties of the class of $N$-harmonic functions is its invariance by conformal transformations of the space  $\BBR^N$. This article is devoted to the study of $N$-harmonic functions which admit an isolated boundary singularity. More precisely, let $a\in\prt\Gw$
 and $u\in W_{loc}^{1,N}(\Gw)\cap C(\overline\Gw\setminus\{a\})$ be a $N$-harmonic function vanishing on $\prt\Gw\setminus\{a\}$, then $u$ may develop a singularity at the point $a$. Our goal is to show the existence of such singular solutions, and then to classify all the positive $N$-harmonic functions with a boundary isolated singularity. We denote by $\bf n_a$ the outward normal unit vector to $\Gw$ at $a$
 The main result we prove are presented below:\\
 
 \noindent {\it There exists a unique positive $N$-harmonic function $u=u_{1,a}$ in $\Gw$, vanishing on $\prt\Gw\setminus\{a\}$ such that
 \begin{equation}\label {behav}
\lim_{\scriptsize\BA{c}x\to a\\
\frac{x-a}{\abs {x-a}}\to\gs\EA}\abs {x-a}u(x)=-\langle\gs,\bf n_a\rangle
\end {equation}
uniformly 
on $S^{N-1}\cap \overline{\Omega}=\{\gs\in S^{N-1}: \langle\gs,{\bf n_{a}}\rangle <0\}$.}\\

The functions $u_{1,a}$ plays a fundamental role in the description of all the positive singular $N$-harmonic functions since we the next result holds\\

 \noindent {\it Let $u$ be a positive $N$-harmonic function  in $\Gw$, vanishing on $\prt\Gw\setminus\{a\}$. Then there exists $k\geq 0$ such that
  \begin{equation}\label {charact}
u=ku_{1,a}.
\end {equation}
}
When $u$ is no longer assumed to be positive we obtain some classification results provided its growth is limited as shows the following\\

 \noindent {\it Let $u$ be a $N$-harmonic function  in $\Gw$, vanishing on $\prt\Gw\setminus\{a\}$ and verifying 
 $$\abs u\leq Mu_{1,a},
 $$ 
 for some $M\geq 0$. Then there exists $k\in\BBR$ such that
  \begin{equation}\label {charact2}
u=ku_{1,a}.
\end {equation}
}

In the last section we give a process to construct $p$-harmonic regular functions ($p>1$) or $N$-harmonic singular functions as product of one variable functions. Starting from the existence of $p$-harmonic functions in the plane under the form $u(x)=u(r,\gs)=r^{\gb}\gw(\gth)$ (see \cite {Kr}), our method, by induction on $N$, allows us to produce separable solutions of the {\it spherical $p$-harmonic spectral equation}
\begin {equation}\label {psrad-p*}
- div_{\gs}\left(\left(\gb^2v^2+\abs {\nabla _{\gs}v}^2\right)^{(p-2)/2}\nabla _{\gs}v\right)
=\gl_{N,\gb}\left(\gb^2v^2+\abs {\nabla _{\gs}v}^2\right)^{(p-2)/2}v.
\end {equation}
on $S^{N-1}$, where $\gl_{N,\gb}=\gb\left(N-1+(\gb-1)(p-1)\right)$. This equation equation is naturally associated to the existence of $p$-harmonic functions under the form $u(x)=\abs x^{\gb}v(x/\abs x)$.
As a consequence,  we express $p$-harmonic functions under the form of a product of $N$-explicit functions of one real variable. {\it If we represent $\BBR^N$ as the set of $\{x=(x_{1},..., x_{N})\}$ where
$x_{1}=r\sin\gth_{N-1}\sin\gth_{N-2}...\sin\gth_{2}\sin\gth_{1}$, 
$x_{2}=r\sin\gth_{N-1}\sin\gth_{N-2}...\sin\gth_{2}\cos\gth_{1}$, ...,
$x_{N-1}=r\sin\gth_{N-1}\cos\gth_{N-2}$ and 
$x_{N}=r\cos\gth_{N-1}$
with $\gth_{1}\in [0,2\gp]$ and $\gth_{k}\in [0,\gp]$, for $k=2,...,N-1$, then, for any integer $k$ the function 
  \begin{equation}\label {charact-p}
u(x)=(r\sin\gth_{N-1}\sin\gth_{N-2}...\sin\gth_{2})^{\gb_{k}}\gw_{k}(\gth_{1})
\end {equation}
is $p$-harmonic in $\BBR^N$, in which expression $\gb_{k}>1$ is an algebraic number depending on $k$ and 
$\gw_{k}$ is a $\gp/k$-antiperiodic solutions of a completely integrable  homogeneous differential equation. Moreover $N$-harmonic singular 
functions are also obtained under the form}
  \begin{equation}\label {charact-p'}
u(x)=r^{-\gb_{k}}(\sin\gth_{N-1}\sin\gth_{N-2}...\sin\gth_{2})^{\gb_{k}}\gw_{k}(\gth_{1}).
\end {equation}

Our paper is organized as follows: 1- Introduction. 2- Construction of fundamental singular $N$-harmonic functions. 3- The classification theorem. 
4- Separable solutions of the $p$-harmonic spectral problem.

\mysection {Construction of fundamental singular $N$-harmonic functions}
We denote by $\CH_N$ the group of conformal transformations in $\BBR^N$. This group is generated by homothethies, inversion and isometries. Our first result is classical, but we repeat the proof for the sake on completeness.
\bprop{inv} Let $u$ be a $N$-harmonic function in a domain $G\subset\BBR^N$
and $h\in \CH_N$. Then $u_h=u\circ h$ is $N$-harmonic in $h^{-1}(G)$.
\es
\Proof Because for any $p>1$ the class of $p$-harmonic functions is invariant by 
homothethies and isometries, it is sufficient to prove the result if $h$ is the inversion 
$\CI_0^1$ with center the origin in $\BBR^N$ and power $1$. We set
$y=\CI_0^1(x)$ and $v(y)=u(x)$. For any $i=1,...,N$
$$u_{x_i}(x)=\sum_j\left(\gd_{ij}\abs x^{-2}-2\abs x^{-4}x_ix_j\right)v_{y_j}(y).
$$
Then
$$\abs{Du}^2(x)=\abs x^{-4}\abs{Dv}^2(y)=\abs y^{4}\abs{Dv}^2(y).
$$
If $\gf$ is a test function, we denote similarly $\psi(y)=\phi(x)$, thus
$$\langle Du,D\phi\rangle=\abs x^{-4}\langle Dv,D\psi\rangle=
\abs y^{4}\langle Dv,D\psi\rangle,
$$
and
$$\int_G\abs {Du}^{N-2}\langle Du,D\phi\rangle\,dx
=\int_{\CI_0^1(G)}\abs y^{2N}\abs {Dv}^{N-2}\langle Dv,D\psi\rangle\abs {D\CI_0^1}\,dy
$$
Because $\abs {D\CI_0^1}=\abs {\det(\prt x_i/\prt y_j)}=\abs y^{-2N}$, the result follows.
\qeda
\\

\bprop {ball} Let $N\geq 2$, $B=B_1(0)$ and $a\in \prt B$. Then there exists a unique positive $N$-harmonic function $U^i$ in $B$ which vanishes on $\prt B\setminus\{a\}$ and satisfies
\begin {equation}\label {behav1}
U^i(x)=\myfrac {1-\abs{x}}{\abs {x-a}^2}(1+\circ (1))\quad \mbox {as }x\to a.
\end {equation}
\es
\Proof We first observe that the coordinates functions $x_i$ are $N$-harmonic and 
positive in the half-space $H_i=\{x\in \BBR^N:x_i>0\}$ and vanishes on 
$\prt H_i$. Therefore, the functions $\chi_i(x)=x_i/\abs x^2$ are also $N$-harmonic and singular at $0$. Without loss of generality we can assume that 
$a$ is the origin of coordinates, 
and that $B$ is the ball with radius $1$ and center $(-1,0,...,0)$. Let 
$\gw$ be the point with coordinates $(-2,0,...,0)$. By
the inversion $\CI_\gw^4$, $a$ is invariant and $B$ is transformed into
the half space $H_1$. Since $\chi_1$ is $N$-harmonic in $H_1$, the function 
$$x\mapsto \chi_1\circ\CI_\gw^4(x)=-\myfrac{\abs x^2+2x_1}{2\abs x^2}$$
is $N$-harmonic and positive in $B=\{x:\abs x^2+2x_1<0\}$, vanishes
on $\prt B$ and is singular at $x=0$. If we set 
$x'_1=x_1+1$,  $x'_i=x_i$ for $i=2,...,N$ and 
$U^i(x')=\chi_1\circ\CI_\gw^4(x)$, then the $x'$ coordinates of $a$ are
$(1,0,...,0)$ and
$$U^i(x')=\myfrac {1-\abs {x'}^2}{2\abs {x'-a}^2}=\myfrac {1-\abs {x'}}{\abs {x'-a}^2}(1+\circ (1))\quad \mbox {as }x'\to a.
$$
Let $\tilde U^i$ be another positive $N$-harmonic function in $B$ which verifies
(\ref {behav1}) and vanishes on $\prt B\setminus\{a\}$. Thus, for any $\gd>0$, $(1+\gd)\tilde U^i$, is positive, $N$-harmonic, and $U^i-(1+\gd)\tilde U^i$ is negative near $a$. By the maximum principle, $U^i\leq (1+\gd)\tilde U^i$. Letting $\gd\to 0$, and permuting $U^i$ and $\tilde U^i$ yields $\tilde U^i=U^i$.\qeda
\\

By performing the inversion $\CI_0^1$, we derive the dual result
\bprop {extball} Let $N\geq 2$, $G=B^c_1(0)$ and $a\in \prt B$. Then there exists a unique positive $N$-harmonic function $U^e$ in $G$ which vanishes on $\prt B\setminus\{a\}$ and satisfies
\begin {equation}\label {behav2}
U^e(x)=\circ (\ln\abs x)\quad \mbox {as }\abs x\to \infty,
\end {equation}
and
\begin {equation}\label {behav3}
U^e(x)=\myfrac {\abs{x}-1}{\abs {x-a}^2}(1+\circ (1))\quad \mbox {as }x\to a.
\end {equation}
\es
\Proof The assumption (\ref{behav2}) implies that the function $U=U^e\circ \CI_0^1$, which is $N$-harmonic in $B\setminus\{0\}$ verifies
$$U(x)=\circ (\ln(1/\abs x))\quad \mbox {near }0.
$$
By \cite {Se1}, $0$ is a removable singularity and thus $U$ can be extended as a positive $N$-harmonic function in $B$ which satisfies 
(\ref{behav1}). This implies the claim.\qeda \\ \hspace{2cm} \\
We denote by $\dot\gr(x)$ the signed distance from $x$ to $\partial \Omega$.
Since $\partial \Omega $ is $C^{2}$, there exists $\beta_{0}>0$ such that if $x\in \mathbb{R}^{N}$ verifies $-\beta_{0} \leq \dot\gr(x)\leq \beta_{0}$, there exists a unique $\xi_{x}\in\prt\Gw$ such that 
 $\abs{x-\xi_{x}}=\abs{\dot\gr(x)}$. Furthermore, if $\gn_{\xi_{x}}$ is the outward unit vector to $\prt\Gw$ at $\xi_{x}$, $x=\xi_{x}-\dot\gr(x)\gn_{\xi_{x}}$. In particular $\xi_{x}-\dot\gr(x)\gn_{\xi_{x}}$ and $\xi_{x}+\dot\gr(x)\gn_{\xi_{x}}$ have the same orthogonal projection $\xi_{x}$ onto $\prt\Gw$. \\ \smallskip
 
 Let $T_{\gb_{0}}(\Gw)=\{x\in\BBR^N:-\gb_{0}\leq\dot\gr(x)\leq\gb_{0}\}$, then the mapping 
 $\Gp:[-\gb_{0},\gb_{0}]\ti\prt\Gw\mapsto T_{\gb_{0}}(\Gw)$ defined by $\Gp(\gr,\xi)=\xi-\gr{\bf\gn}(\xi)$ is a $C^2$ diffeomorphism. Moreover $D\Gp(0,\xi)(1,e)=e-\gn_{\xi}$ for any $e$ belonging to the tangent space $T_{\xi}(\prt\Gw)$ to $\prt\Gw$ at $\xi$. If $x\in T_{\gb_{0}}(\Gw)$, we define the reflection of $x$ through $\prt\Gw$ by
 $\psi(x)=\xi_{x}+\dot\gr(x)\gn_{\xi_{x}}$. Clearly $\psi$ is an involutive diffeomorphism from $\overline\Gw\cap T_{\gb_{0}}(\Gw)$ to $\Gw^c\cap T_{\gb_{0}}(\Gw)$, and $D\psi (x)=I$ for any $x\in\prt\Gw$.
If a function $v$ is defined in $\Gw\cap T_{\gb_{0}}(\Gw)$, we define $\tilde v$ in $ T_{\gb_{0}}(\Gw)$ by\\
 \begin {equation}\label {ext}
\tilde v(x)=\left\{\BA{l}v(x)\quad\qquad\mbox {if }x\in \Gw\cap T_{\gb_{0}}(\Gw)\\[2mm]
-v\circ\psi(x)\quad\mbox {if }x\in \Gw^c\cap T_{\gb_{0}}(\Gw).
\EA\right.\end {equation}
\\
\blemma{ext} Assume that $0\in \partial \Omega$. Let $v\in C^{1,\ga}(\overline\Gw\cap T_{\gb_{0}}(\Gw)\setminus\{0\})$ be a solution of (\ref {p-harm}) in $\Gw\cap T_{\gb_{0}}(\Gw)$ vanishing on $\prt\Gw\setminus\{0\}$. Then 
$\tilde v\in C^{1,\ga}(T_{\gb}(\Gw)\setminus\{0\})$ is solution of a quasilinear equation
 \begin {equation}\label {ext-equ}
\mysum{j}{}\myfrac {\prt }{\prt x_{j}}\tilde A_{j}(x,D\tilde v)=0
\end {equation}
in $T_{\gb}(\Gw)\setminus\{0\}$ where the $\tilde A_{j}$ are $C^1$ functions defined in $T_{\gb}(\Gw)$ where they verify 
 \begin {equation}\label {ext1}\left\{\BA {l}
 (i)\quad \tilde A_{j}(x,0)=0\\[2mm]
 
 (ii) \quad\mysum{i,j}{}\myfrac {\prt \tilde A_{j}}{\prt \eta_{i}}(x,\eta)\xi_{i}\xi_{j}
 \geq \Gg\abs\eta^{p-2}\abs \xi^2
 \\
 
 (iii) \quad\mysum{i,j}{}\abs {\myfrac {\prt  \tilde A_{j}}{\prt \eta_{i}}(x,\eta)}\leq \Gg\abs\eta^{p-2}\\
 
\EA\right.
\end {equation}
for all $x\in T_{\gb}(\Gw)\setminus\{0\}$ for some $\gb\in (0,\gb_{0}]$, $\eta\in \BBR^N$, $\xi\in\BBR^N$ and some 
$ \Gg>0$.
\es
\Proof The assumptions (\ref {ext1}) implies that weak solutions of (\ref {ext-equ}) are $C^{1,\ga}$, for some $\ga>0$ \cite {To1} and satisfy the standard a priori estimates. As it is defined the function $\tilde v$ is clearly $C^1$ in $T_{\gb_{0}}(\Gw)\setminus\{0\}$. Writing
$Dv(x)=-D(\tilde v\circ \psi(x))=-D\psi(x)(D\tilde v(\psi(x)))$ and $\tilde x=\psi(x)=\psi^{-1}(x)$
$$\BA {l}
\myint {\Gw\cap T_{\gb}(\Gw)}{}\abs {Dv}^{p-2}Dv.D\gz dx\\

\phantom {-----}
=\myint {\overline\Gw^c\cap T_{\gb}(\Gw)}{}\abs {D\psi(D\tilde v)}^{p-2}
D\psi (D\tilde v).D\psi(D\gz)\abs {D\psi}d\tilde x.
\EA$$
But 
$$\BA {l}D\psi (D\tilde v).D\psi(D\gz)=\mysum{k}{}
\left(\mysum{i}{}\myfrac {\prt \psi_{i}}{\prt x_{k}}\myfrac {\prt \tilde v}{\prt x_{i}}\right)
\left(\mysum{j}{}\myfrac {\prt \psi_{j}}{\prt x_{k}}\myfrac {\prt \gz}{\prt x_{j}}\right)\\
\phantom{D\psi (D\tilde v).D\psi(D\gz)}
=\mysum{j}{}\left(\mysum{i,k}{}\myfrac {\prt \psi_{i}}{\prt x_{k}}\myfrac {\prt \psi_{j}}{\prt x_{k}}\myfrac {\prt \tilde v}{\prt x_{i}}\right)\myfrac {\prt \gz}{\prt x_{j}}.
\EA$$
We set $b(x)=\abs{D\psi}$, 
\begin {equation}\label {ext2}
A_j(x,\eta)=\abs {D\psi}\abs {D\psi(\eta)}^{p-2}\mysum{i}{}\left(\mysum{k}{}\myfrac {\prt \psi_{i}}{\prt x_{k}}\myfrac {\prt \psi_{j}}{\prt x_{k}}\right)\eta_i,
\end {equation}
and
\begin {equation}\label {ext3}
A(x,\eta)=(A_1(x,\eta),...,A_N(x,\eta))=\abs {D\psi}\abs {D\psi(\eta)}^{p-2}(D\psi)^tD\psi(\eta).
\end {equation}
For any $\xi\in\prt\Gw$, the mapping
$D\psi_{\prt\Gw}(\xi)$ is the symmetry with respect to the hyperplane $T_{\xi}(\prt\Gw)$ tangent to $\prt\Gw$ at $\xi$, so $\abs {D\psi(\xi)}=1$. Inasmuch $D\psi$ is continuous, a lengthy but standard computation leads to the existence of some $\gb\in (0,\gb_{0}]$ such that (\ref {ext1}) holds in $T_{\gb}(\Gw)\cap\overline\Gw^c$. If we define $\tilde A$ to be $\abs\eta^{p-2}\eta$ on $T_{\gb}(\Gw)\cap\overline\Gw$ and $A$  on $T_{\gb}(\Gw)\cap\overline\Gw^c$, then inequalities (\ref {ext1}) are satisfied in $T_{\gb}(\Gw)$.\qeda 
\\

These three results allows us to prove our main result
\bth{gen} Let $\Gw$ be an open subset of $\BBR^N$ with a compact $C^2$ boundary, $\gr(x)=\dist (x,\prt\Gw)$ and $a\in\prt \Gw$. Then there exists one and only one positive $N$-harmonic function $u$ in $\Gw$, vanishing on 
$\prt\Gw\setminus\{a\}$ verifying
\begin {equation}\label {behav4}
\lim_{\scriptsize\BA{c}x\to a\\
\frac{x-a}{\abs {x-a}}\to\gs\EA}\abs {x-a}u(x)=-\langle\gs,\bf n_a\rangle
\end {equation}
uniformly on $S^{N-1}\cap \overline{\Omega}$, and 
\begin {equation}\label {behav5}
u(x)=\circ (\ln\abs x))\quad \mbox {as }\abs x\to \infty,
\end {equation}
if $\Gw$ is not bounded.
\es
\Proof Uniqueness follows from (\ref{behav4}) by the same technique as in the previous propositions.\\
{\it Step 1 (Existence). }If $\Gw$ is not bounded, we perform an inversion 
$\CI_m^{\abs{m-a}^2}$ with center some $m\in\Gw$. Because of (\ref{behav5}), the new function $u\circ \CI_m^{\abs{m-a}^2}$ is $N$-harmonic in $\Gw'=\CI_m^{\abs{m-a}^2}(\Gw)$ and satisfies (\ref{behav4}). Thus we are reduced to the case were $\Gw$ is bounded. Since $\Gw$ is $C^2$, it satisfies the interior and exterior sphere condition at $a$. By dilating $\Gw$, we can assume that the exterior and interior tangent spheres at $a$ have radius $1$. We denote them by $B_1(\gw^e)$ and $B_1(\gw^i)$, their respective centers being
$\gw^i=a-\bf n_a$ and $\gw^e=a+\bf n_a$. We set
$V^i(x)=U^i(x-\gw^i)$ and $V^e(x)=U^e(x-\gw^e)$ where $U^i$ and 
$U^e$ are the two singular $N$-harmonic functions described in 
\rprop {ball} and \rprop {extball}, respectively in $B_1(\gw^i)$ and 
$B^c_1(\gw^e)$, with singularity at point $a$. For 
$\ge>0$, we put $\Gw_\ge=\Gw\setminus B_\ge(a)$, 
$\Gs_\ge=\Gw\cap\prt B_\ge(a)$ and 
$\prt^* \Gw_\ge=\prt \Gw\cap B^c_\ge(a)$. Let $u_\ge$ be the solution of
\begin {equation}\label {approx}\left\{\BA{l}
{div}(\abs {Du_\ge}^{N-2}Du_\ge)=0\quad\mbox {in }\Gw_\ge\\
\phantom {-------;,}
u_\ge=0\quad\mbox {on }\prt^* \Gw_\ge\\
\phantom {-------;,}
u_\ge=V^e\quad\mbox {on }\Gs_\ge.
\EA\right.\end {equation}
This solution is obtained classicaly by minimisation of a convex functional over a class of functions with prescribed boudary value on
$\prt\Gw_\ge$. For any $x\in B_1(\gw^i)$, there holds
$$\dist (x,\prt B_1(\gw^e))=\abs {x-\gw^e}-1\geq 
\dist (x,\prt \Gw)\geq \dist (x,\prt B_1(\gw^i))=1-\abs {x-\gw^i}.
$$
thus
$$V^i(x)\leq V^e(x)\forevery x\in B_1(\gw^i),
$$
by using (\ref {behav1}), (\ref {behav3}) and the maximum principle. Therefore
$$V^i(x)\leq u_\ge (x)\leq V^e(x)\forevery x\in B_1(\gw^i)\cap \Gw_\ge
$$
and
$$u_\ge (x)\leq V^e(x)\forevery x\in \Gw_\ge.
$$
Finally, for $0<\ge'<\ge$, 
$u_{\ge'}\vline_{\Gs_\ge}\leq V^e \vline_{\Gs_\ge}= u_{\ge}\vline_{\Gs_\ge}$. Thus 
$$u_{\ge'}(x)\leq u_{\ge}(x)\forevery x\in\Gw_\ge.
$$
The sequence $\{u_\ge\}$ is increasing with $\ge$. By classical a priori estimates concerning quasilinear equations, it converges to some positive $N$-harmonic function $u$ in $\Gw$ which vanishes on 
$\prt\Gw\setminus\{a\}$ and verifies
$$V^i(x)\leq u (x)\forevery x\in  B_1(\gw^i),
$$
and
$$u (x)\leq U^e(x)\forevery x\in  \Gw.
$$
This implies
\begin {equation}\label {approx1}
\myfrac {1-\abs {x-\gw_i}^2}{2\abs {x-a}^2}\leq u (x)\forevery x\in  B_1(\gw^i),
\end {equation}

\begin {equation}\label {approx2}
 u (x)\leq\myfrac {\abs {x-\gw_e}^2-1}{2\abs {x-a}^2}\forevery x\in  \Gw,
\end {equation}
By scaling we can prove the following estimate
\begin {equation}\label {approx3}
u (x)\leq C\myfrac {\gr(x)}{\abs {x-a}^2}\forevery x\in  \Gw.
\end {equation}
for some $C>0$: for simplicity we can assume that $a$ is the origin of coordinates and, for $r>0$ set
$u_r(y)=u(r y)$. Clearly $u_r$ is $N$-harmonic in $\Gw/r$ and 
$$\max \{\abs {Du_r(y)}:y\in \Gw/r\cap (B_{3/2}\setminus B_{2/3}) \}
\leq C\max \{\abs {u_r(z)}:z\in \Gw/r\cap (B_2\setminus B_{1/2})\},
$$
where $C$, which depends on the curvature of $\prt\Gw/r$, remains bounded as long as $r\leq 1$.
Since $Du_r(y)=r Du(r y)$, we obtain by taking $r y=x$, $\abs y=1$ and using (\ref{approx2}) with general $a$, $ \abs {Du(x)}\leq C\abs {x-a}^{-2}$. By the mean value theorem, since $u$ vanishes on $\prt\Gw\setminus \{a\}$, (\ref{approx3}) holds.\\
{\it  Step 2. } In order to give a simple proof of the estimate (\ref{behav4}), we fix the origin of coordinates at $a=0$ and the normal outward unit vector at $a$ to be $-\mathbf{e}_{N}$. If $\tilde{u}$ is the extension of $u$ by reflection through $\prt \Omega$, it statisfies (\ref{ext-equ}) in $T_{\gb}(\Gw)\backslash\{0\}$ (see lemma \ref{ext}). For $r>0$, set $\tilde{u}^{r}(x)=r \tilde{u}(rx)$. Then $\tilde{u}^{r}$ is solution of
\begin{equation}\label{eq}
\mysum{j}{}\myfrac {\prt }{\prt x_{j}}\tilde A_{j}(rx,D\tilde u^{r})=0
\end{equation}
in $T_{\gb/r}(\Gw/r)\backslash\{0\}$. By the construction of $\tilde{A}_{j}(x,\eta)$, we can note that\\
\begin{center}
$\displaystyle{\lim_{\scriptsize r\to 0} \tilde{A}^{j}(rx,\eta)=|\eta|^{p-2}\eta_{j},\hspace{1cm}\forall \eta \in \mathbb{R}^{N}}$.
\end{center}
Furthermore, for any $x\in T_{\beta}(\Gw)\backslash\{0\}$, $\rho(x)=\rho(\psi(x))$ and $c\abs{x}\leq \abs{\psi(x)}\leq c^{-1} \abs{x}$ for some $c>0$, the estimate (\ref{approx3}) holds if $u$ is replaced by $\tilde{u}^{r}$, $\Gw $ by $T_{\gb/r}(\Gw/r)$ and $\rho(x)$ by $\rho_{r}(x):=$ dist$(x,\Gw/r)$ i.e.
\begin{center}
$\abs{\tilde{u}^{r}(x)}\leq C |x|^{-2} \rho_{r}(x)$ $\forall x\in T_{\gb/r}(\Gw/r)$.
\end{center}
For $0<a<b$ fixed and for some $0< r_{0}\leq 1$ the spherical shall $\Gamma_{a,b}=\{x\in \mathbb{R}^{N} : a\leq \abs{x}\leq b \}$ is included into $T_{\beta/r}(\Gw/r)$ for all $0<r\leq r_{0}$. By the classical regularity theory for quasilinear equations \cite{To1} and lemma \ref{ext}, there holds
\begin{center}
$\left\| D\tilde{u}^{r} \right\|_{C^{\alpha}(\Gamma_{2/3,3/2})}\leq C_{r} \left\|\tilde{u}^{r}\right\|_{L^{\infty}(\Gamma_{1/2,2})},$
\end{center}
where $C_{r}$ remains bounded because $r\leq 1$. By Ascoli's theorem, (\ref{approx1}) and (\ref{approx3}), $\tilde{u}^{r}(x)$ converges to $x_{N}|x|^{-2}$ in the $C^{1}(\Gamma_{2/3,3/2})$-topology. This implies in particular that $ r^{2}D\tilde{u}(rx)$ converges uniformly in $\Gamma_{2/3,3/2}$ to $-2 x_{N}|x|^{-4}x+|x|^{-2} \mathbf{e}_{N}$. Using the expression of $D\tilde{u}$ in spherical coordinates we obtain
\begin{center}
$r^{2}\tilde{u}_{r} \mathbf{i}- r \tilde{u}_{\phi} \mathbf{e}+ \myfrac{r}{\sin \phi}\nabla_{\sigma'}\tilde{u} \to\ -2 \sigma_{N} \mathbf{i} + \mathbf{e}_{N}$ uniformly on $S^{N-1}$ as $r\to\ 0$,
\end{center}
where $\cos \phi = x_{N}|x|^{-1}$, $\mathbf{i}=x/|x|$, $\mathbf{e}$ is derived from $x/|x|$ by a rotation with angle $\pi/2$ in the plane $0,x,N$ ($N$ being the North pole), and $\nabla_{\sigma'} $ is the covariant gradient on $S^{N-2}$. Inasmuch $\mathbf{i}$, $\mathbf{e}$ and $\nabla_{\sigma'}$ are orthogonal, the components of $\mathbf{e}_{N}$ are $\cos \phi, \sin \phi$ and $0$, thus
\begin{center}
$r \tilde{u}_{\phi}(r,\sigma',\phi) \to\ - \sin \phi $ as $r \to\ 0$.
\end{center}
Since
\begin{center}
$\tilde{u}(r,\sigma',\phi)=\myint{\pi/2}{\phi}\tilde{u}_{\phi}(r,\sigma',\theta)d\theta$,
\end{center}
the previous convergence estimate establishes (\ref{behav4}).\qeda \\

\bdef{fund} We shall denote by $u_{1,a}$ the unique positive $N$-harmonic function satisfying (\ref{behav4}), and call it the fundamental solution with a point singularity at $a$.
\es

\mysection {The classification theorem}
In this section we characterize all the positive $N$-harmonic functions
vanishing on the boundary of a domain except one point. The next statement is an immediate consequence of \rth {gen} and \cite [Th. 2.11]{BBV}.
\bth {est}. Let $\Gw$ be a bounded domain with a $C^2$ boundary and $a\in\prt\Gw$. If $u$ is a positive $N$-harmonic function in $\Gw$ vanishing on $\prt\Gw\setminus\{a\}$, there exists $M\geq 0$ such that 
\begin {equation}\label {est1}
u(x)\leq Mu_{1,a}(x)\forevery x\in\Gw
\end {equation}
\es

In the next theorem, which extends \cite [Th. 2.13]{BBV}, we characterize all the signed $N$-harmonic functions with a moderate growth near the singular point. 

\bth {charact}. Let $\Gw$ be a bounded domain with a $C^2$ boundary and $a\in\prt\Gw$. Assume that $u_{1,a}$ has only a finite number of critical points in $\Gw$. If $u$ is a $N$-harmonic function in $\Gw$ vanishing on $\prt\Gw\setminus\{a\}$ verifying 
$\abs {u(x)}\leq Mu_{1,a}(x)$ for some $M>0$ and any $x\in\Gw$, there
exists $k\in[-M,M]$ such that $u=ku_{1,a}$.
\es
\Proof We define $k$ as the minimum of the $\ell$ such that $u\leq \ell u_{1,a}$ in $\Gw$. Without any loss of generality we can assume $k>0$. Then either the tangency of the graphs of the functions $u$ and $ku_{1,a}$ is achieved in $\overline\Gw\setminus \{a\}$, or it is achieved asymptotically at the singular point $a$. In the first case we considered two sub-cases:\\[1mm]
(i) The coincidence set $G$ of $u$ and $ku_{1,a}$ has a connected component $\gw$ isolated in $\Gw$. In this case there exists a smooth domain $\CU$ such that $\overline\gw\subset\CU$ and $\gd>0$ such that $ku_{1,a}-u\geq\gd$ on $\prt\CU$. The maximum principle implies that $ku_{1,a}-u\geq\gd$ in $\CU$, a contradiction.\\[1mm]
(ii) In the second sub-case any connected component $\gw$ of the coincidence set touches $\prt\Gw\setminus\{a\}$, or the two graphs admits a tangency point on
$\prt\Gw\setminus\{a\}$. If $m\in\gw\cap\prt\Gw\setminus\{a\}$ or is such a tangency point, the regularity theory implies $\prt u(m)/\prt{\bf n}_m=ku_{1,a}(m)/\prt{\bf n}_m$. By Hopf boundary lemma,  $u_{1,a}(m)/\prt{\bf n}_m<0$. By the mean value theorem, the function $w=ku_{1,a}-u$ satisfies an equation
\begin {equation}\label {lin}
Lw=0
\end {equation}
which is elliptic and non degenerate near $m$ (see \cite {FV}, \cite {KV}), it follows that
$w$ vanishes in a neighborhood of $m$ and the two graphs cannot be tangent only on 
$\prt\Gw\setminus\{a\}$. Assuming that $\gw\neq\Gw$, let $x_0\in\Gw\setminus\gw$
such that $\dist (x_0,\gw)=r_0<\gr(x_0)=\dist (x_0,\prt\Gw$, and let $y_0\in\gw$ be such that $\abs {x_0-y_0}=r_0$. Since $u_{1,a}$ has at most a finite number of critical points, we can choose $x_0$ such that $y_0$ is not one of these critical points. By assumption $w=ku_{1,a}-u$ is positive in $B_{r_0}(x_0)$ and vanishes at  a boundary point $y_0$. Since the equations are not degenerate at $y_0$ there holds
$$k\prt u_{1,a}(y_0)/\prt {\bf \gn}-\prt u(y_0)/\prt {\bf \gn}< 0
$$
where $\gn=(y_0-x_0)/r_0$, which contradicts the fact that the two graphs are tangent at $y_0$.\smallskip

\noindent Next we are reduced to the case where the graphs of $u$ and $ku_{1,a}$ are separated in $\Gw$ and asymptotically tangent at the singular point $a$. There exists a sequence $\{\xi_n\}\subset \Gw$ such that 
$\lim_{n\to\infty} u(\xi_n)/u_{1,a}(\xi_n)=k$. We set 
$\abs {x_{n}-a}=r_{n}$, $u_{n}(y)=r_{n}u(a+r_{n}y)$ and $v_{n}(y)=r_{n}u_{1,a}(a+r_{n}y)$. Both $u_{n}$ and $ v_{n}$ are $N$-harmonic in $\Gw_{n}=(\Gw-a)/r_{n}$. The functions $u_{n}$ and $v_{n}$ are locally uniformly bounded in $\overline\Gw_{n}\setminus \{0\}$. It follows, by using classical regularity results, that,  there exists  sub-sequences, such that  $\{u_{n_k}\}$ and $\{v_{n_k}\}$ converge respectively to $U$ and $V$
in the $C^1_{loc}$-topology of $\overline\Gw_{n_{k}}\setminus\{0\}$. The functions $U$ and $V$ are $N$-harmonic in $H\approx\BBR^N_+=\{x=(x_1,x_2,...,x_N):x_N>0\}$ and vanish on $\prt H\setminus\{0\}$. Since it can be assumed that $(\xi_{n_k}-a)/r_{n_k}\to \xi$, there holds
$U\leq kV$ in $H$, $U(\xi)=kV(\xi)$, if $\xi\in H$, and $\prt U(\xi)/\prt x_N=k\prt V(\xi)/\prt x_N>0$, if 
$\xi\in \prt H$ (notice that $\abs {\xi}=1$). If $\xi\in\prt H$, Hopf lemma applies to $V$ at $\xi$ and, using the same linearization with the linear operator $L$ as in the previous proof, it yields to
$U=kV$. If $\xi\in H$, we use the fact that $\abs {Du_{1,a}(x)}\geq \gb>0$ for $\abs {x-a}\leq \ga$ for some $\gb,\ga>0$. Thus $\abs {Dv_n(\xi)}\geq \gb$. The non-degeneracy of $V$ and the strong maximum principle lead again to $U=kV$. Whatever is the position of $\xi$, the equality between $U$ and $kV$ and the convergence in 
$C^1_{loc}$ leads to the fact that for any $\ge>0$ there exists $n_\ge\in\BBN$ such that $n\geq n_\ge$ implies 
$$(k-\ge)u_{1,a}(x)\leq u(x)\leq (k+\ge)u_{1,a}(x)\forevery x\in\Gw\cap \prt B_{r_n}(a).
$$
By the comparison principle between $N$-harmonic functions this inequality holds true in $\Gw\setminus \prt B_{r_n}(a)$. Since $r_n\to 0$ and $\ge$ is arbitrary, this ends the proof. 
\qeda \\

\noindent\Remark The assumption that $u_{1,a}$ has only isolated critical points in 
$\Gw$ is clearly satisfied in the case of a ball, a half-space or the complementary of a ball where no critical point exists. It is likely that this assumption always holds but we cannot prove it. However the Hopf maximum principle for p-harmonic functions (see \cite {To}) implies that $u_{1,a}$ cannot have local extremum in $\Gw$.
\mysection {Separable solutions of the $p$-harmonic spectral problem}
In this section we present  a technique for constructing  signed $N$-harmonic functions, regular or singular,
as a product of  functions depending only on one real variable. Some of the results were sketched in \cite {Ve4}. The starting point is the result of Krol \cite {Kr} dealing with the existence of $2$-dimensional separable $p$-harmonic functions (the construction of singular separable $p$-harmonic functions was performed in \cite {KV}). 
\bth {Krth} {\rm (Krol)} Let  $p>1$. For any positive integer $k$ there exists a unique $\gb_{k}>0$ and 
$\gw_{k}:\BBR\mapsto\BBR$, with least antiperiod $\gp/k$, of class $C^\infty $ such that
\begin {equation}\label {psrad1}
u_{k}(x)=\abs x^{\gb_{k}}\gw_{k}(x/\abs x)
\end {equation}
is $p$-harmonic in $\BBR^2$; $\gb_{k}$  is the unique root $\geq 1$ of
\begin {equation}\label {psrad2}
(2k-1)X^2-\myfrac {pk^2+(p-2)(2k-1)}{p-1}X+k^2=0.
\end {equation}
$(\gb_{k},\gw_{k})$ is unique up to translation and homothety over $\gw_{k}$.
\es
This result is obtained by solving the homogeneous differential equation satisfied by $\gw_{k}=\gw$:
\begin {equation}\label {psrad-1}
-\left(\left(\gb^2\gw^2+\gw^2_{\gth}\right)^{(p-2)/2}\gw_{\gth}\right)_{\gth}
=\gb\left(1+(\gb-1)(p-1)\right)\left(\gb^2\gw^2+\gw^2_{\gth}\right)^{(p-2)/2}\gw.
\end {equation}
In the particular case $k=1$, then $\gb_{1}=1$ and $\gw_{1}(\gth)=\sin\gth$.  For the other values of $k$ the $\gb_{k}$ are algebraic numbers and the $\gw_{k}$ are not trigonometric functions, except if $p=2$. More generally, if one looks for $p$-harmonic functions in $\BBR^N\setminus\{0\}$ under the form
$u(x)=u(r,\gs)=r^\gb v(\gs)$, $r=\abs x>0$, $\gs=x/\abs x\in S^{N-1}$, one obtains that $v$ verifies
\begin {equation}\label {psrad-p}
- div_{\gs}\left(\left(\gb^2v^2+\abs {\nabla _{\gs}v}^2\right)^{(p-2)/2}\nabla _{\gs}v\right)
=\gl_{N,\gb}\left(\gb^2v^2+\abs {\nabla _{\gs}v}^2\right)^{(p-2)/2}v
\end {equation}
on $S^{N-1}$, where $\gl_{N,\gb}=\gb\left(N-1+(\gb-1)(p-1)\right)$ and $div_{\gs}$ and $\nabla_{\gs}$ are respectively the divergence and the gradient operators on $S^{N-1}$ (endowed with the Riemaniann structure induced by the imbedding of the sphere into $\BBR^N$). This equation, called the {\it spherical $p$-harmonic spectral problem}, is the natural generalization of the spectral problem of the Laplace-Beltrami operator on $S^{N-1}$. Since it does not correspond to a variational form (except if $p=2$), it is  difficult to obtain solutions. In the range of $1<p\leq N-1$, Krol proved in \cite {Kr} the existence of solutions of (\ref {psrad-p}), not on the whole sphere, but on a spherical cap (which reduced (\ref {psrad-p}) to an non-autonomous nonlinear second order differential equation). His methods combined ODE estimates and shooting arguments. Later on, Tolksdorf \cite {To} introduced an entirely new method for proving the existence of solutions on any $C^2$ spherical domain $S$, with Dirichlet boundary conditions. Only the case $\gb>0$ was treated in \cite {To}, and, by a small adaptation of Tolksdorf approach, the case $\gb>0$ was considered in \cite {Ve4}. We develop below a method which allows to express solutions as product of explicit one variable functions.

\subsection {The $3$-D case}
Let $(r,\gth,\gf)\in (0,\infty)\ti[0,2\gp]\ti [0,\gp]$ be the spherical coordinates in $\BBR^3$
$$\left\{\BA {l}
x_{1}=r\cos\gth\sin\gf\\
x_{2}=r\sin\gth\sin\gf\\
x_{3}=r\cos\gf
\EA\right.$$
Then (\ref {psrad-p}) turns into
\begin {equation}\label {psrad-p-3}\BA {c}
- \myfrac {1}{\sin\gf}\left[\sin\gf\left(\gb^2v^2+v_{\gf}^2+\myfrac {v^2_{\gth}}{\sin^2\gf}\right)^{(p-2)/2}\!\!\!\!\!\!\!\!\!\!\!\!\!\!\!\!\!v_{\gf}\;\;\;\;\;\;\right]_{\gf}-\myfrac {1}{\sin^2\gf}
\left[\left(\gb^2v^2+v_{\gf}^2+\myfrac {v^2_{\gth}}{\sin^2\gf}\right)^{(p-2)/2}
\!\!\!\!\!\!\!\!\!\!\!\!\!\!\!\!\!v_{\gth}\;\;\;\;\;\;\right]_{\gth}\\[6mm]
\phantom {--------}
=\gb\left(2+(\gb-1)(p-1)\right)\left(\gb^2v^2+v_{\gf}^2+\myfrac {v^2_{\gth}}{\sin^2\gf}\right)^{(p-2)/2}\!\!\!\!\!\!\!\!\!\!\!\!\!\!\!\!\!v
\EA\end {equation}
We look for a function $v$ under the form
\begin {equation}\label {psrad-4}
v(\gth,\gf)=(\sin\gf)^\gb \gw(\gth)
\end {equation}
then
$$\gb^2v^2+v_{\gf}^2+\myfrac {v^2_{\gth}}{\sin^2\gf}
=(\sin\gf)^{2\gb-2}(\gb^2 \gw^2+\gw_{\gth}^2),
$$

$$\myfrac {1}{\sin^2\gf}
\left[\left(\gb^2v^2+v_{\gf}^2+\myfrac {v^2_{\gth}}{\sin^2\gf}\right)^{(p-2)/2}
\!\!\!\!\!\!\!\!\!\!\!\!\!\!\!\!\!v_{\gth}\;\;\;\;\;\;\right]_{\gth}
=(\sin\gf)^{(\gb-1)(p-1)-1}\left((\gb^2 \gw^2+\gw_{\gth}^2)^{(p-2)/2}\gw_{\gth}\right)_{\gth},
$$

$$\left(\gb^2v^2+v_{\gf}^2+\myfrac {v^2_{\gth}}{\sin^2\gf}\right)^{(p-2)/2}\!\!\!\!\!\!\!\!\!\!\!\!\!\!\!\!\!v
=(\sin\gf)^{(\gb-1)(p-1)+1}(\gb^2 \gw^2+\gw_{\gth}^2)^{(p-2)/2}\gw,
$$
and
$$\BA {l}
 \myfrac {1}{\sin\gf}\left[\sin\gf\left(\gb^2v^2+v_{\gf}^2+\myfrac {v^2_{\gth}}{\sin^2\gf}\right)^{(p-2)/2}\!\!\!\!\!\!\!\!\!\!\!\!\!\!\!\!\!v_{\gf}\;\;\;\;\;\;\right]_{\gf}\\[6mm]

 =\gb(\sin\gf)^{(\gb-1)(p-1)-1}
 \left[((\gb-1)(p-1)+1)-\sin^2\gf\;((\gb-1)(p-1)+2)\right](\gb^2 \gw^2+\gw_{\gth}^2)^{(p-2)/2}\gw.
 \EA
 $$
It follows that  $\gw$ satisfies the same equation (\ref {psrad-1}). The next result follows immediately from \rth {Krth}
\bth {septh} Assume $N=3$ and $p>1$. Then for any positive integer $k$ there exists a $p$-harmonic function $u$ in $\BBR^3$ under the form
\begin {equation}\label {psrad-5}
u(x)=u(r,\gth,\gf)=r^{\gb_{k}}(\sin\gf)^{\gb_{k}}\gw_{k}(\gth)
\end {equation}
where $\gb_{k}$ and $\gw_{k}$ are as in \rth {Krth}.
\es
In the case $p=3$ we can use the conformal invariance of the $3$-harmonic equation in $\BBR^3$ to derive
\bth {septh} Assume $p=N=3$. Then for any positive integer $k$ there exists a $p$-harmonic function $u$ in $\BBR^3\setminus\{0\}$ under the form
\begin {equation}\label {psrad-5'}
u(x)=u(r,\gth,\gf)=r^{-\gb_{k}}(\sin\gf)^{\gb_{k}}\gw_{k}(\gth)
\end {equation}
where $\gb_{k}$ and $\gw_{k}$ are as in \rth {Krth} with $p=3$.
\es

As a consequence of \rth {septh} we obtain signed $3$-harmonic functions under the form
(\ref {psrad-5}) in the half space $\BBR^3_{+}=\{x:x_2>0\}$, vanishing on $\prt\BBR^3_{+}\setminus\{0\}$, with a singularity at $x=0$. They correspond to even integers $k$. The extension to general smooth domains $\Gw$ is a deep chalenge. In the particular case $k=1$, we have seen that $\gb_{1}=1$ and $\gw_{1}(\gth)=\sin\gth=x_{2}$, that we already know.
\subsection {The general case}
We assume that $N>3$ and write the spherical coordinates in $\BBR^N$ under the form
\begin {equation}\label {repres0}
x=\left\{(r,\gs)\in (0,\infty)\ti S^{N-1}=(r,\sin\gf\,\gs',\,\cos\gf):\gs'\in S^{N-2},\gf\in [0,\gp]\right\}.
\end {equation}
The main result concerning separable $p$-harmonic functions is the following.
\bth {repres} Let $N>3$ and $p>1$. For any positive integer $k$ there exists $p$-harmonic functions
in $\BBR^N$ under the form
\begin{equation}\label{psrad-N}
u(x)=u(r,\gs',\gf)=(r\sin\gf)^{\gb_{k}}\, w(\gs').
\end{equation}
where $\gb_{k}$  is the unique root $\geq 1$ of (\ref {psrad2}) and $w$ is solution of (\ref {psrad-N-1})
with $\gb=\gb_{k}$. Furthermore, if $p=N$ there exists a singular $N$-harmonic function under the form
\begin{equation}\label{psrad-N=p}
u(x)=u(r,\gs',\gf)=r^{-\gb_{k}}(\sin\gf)^{\gb_{k}}\, w(\gs').
\end{equation}
\es
\Proof We first recall (see \cite {Vi} for details) that the $SO(N)$ invariant unit measure on $S^{N-1}$ is $d\gs=a_{N}
\sin^{N-2}\gf\,d\gs'$ for some $a_{N}>0$, and
$$\nabla_{\gs}v=-v_{\gf}{\bf e}+\myfrac {1}{\sin\gf}\nabla_{\gs'}v.
$$
where $\bf e$ is derived from $x/\abs x$ by the rotation of center $0$ angle $\gp/2$ in the plane going thru $0$, $x/\abs x$ and the north pole. The weak formulation of (\ref{psrad-p}) expresses as
\begin {equation}\label {psrad-p2}\BA{l}
\myint{0}{\gp}\myint{S^{N-2}}{}\!\!
\left(\gb^2v^2+v_{\gf}^2+\myfrac {1}{\sin^2\gf}\abs{\nabla_{\gs'}v}^2\right)^{(p-2)/2}\!\!\left(v_{\gf}\gz_{\gf}+
\myfrac {1}{\sin^2\gf}\nabla_{\gs'}v.\nabla_{\gs'}\gz
\right)\sin^{N-2}\gf\,d\gs'\,d\gf\\[4mm]
\phantom{-------}
=\gl_{N,\gb}\myint{0}{\gp}\myint{S^{N-2}}{}\!\!
\left(\gb^2v^2+v_{\gf}^2+\myfrac {1}{\sin^2\gf}\abs{\nabla_{\gs'}v}^2\right)^{(p-2)/2}\!\!v\,\gz\sin^{N-2}\gf\,d\gs'\,d\gf
\EA\end {equation}
or, equivalently
\begin {equation}\label {psrad-p3}\BA{l}
-\myfrac {1}{\sin^{N-2}\gf}\left[\sin^{N-2}\gf\left(\gb^2v^2+v_{\gf}^2+\myfrac {1}{\sin^2\gf}\abs{\nabla_{\gs'}v}^2\right)^{(p-2)/2}\!\!v_{\gf}\right]_{\gf}
\\[6mm]
\phantom{---}
-\myfrac {1}{\sin^2\gf}div_{\gs'}\left[\left(\gb^2v^2+v_{\gf}^2+\myfrac {1}{\sin^2\gf}\abs{\nabla_{\gs'}v}^2\right)^{(p-2)/2}\!\!\nabla_{\gs'}v\right]
\\[6mm]
\phantom{----------}
=\gl_{N,\gb}
\left(\gb^2v^2+v_{\gf}^2+\myfrac {1}{\sin^2\gf}\abs{\nabla_{\gs'}v}^2\right)^{(p-2)/2}\!\!v
\EA\end {equation}
where $div_{\gs'}$ is the divergence operator acting on vector fields on $S^{N-2}$. We look again for p-harmonic functions under the form
\begin{equation}\label{psrad-N0}
u(r,\gs)=u(r,\gs',\gf)=r^\gb v(\gs',\gf)=r^\gb\sin^\gb\gf\, w(\gs').
\end{equation}
Then
$$\left(\gb^2v^2+v_{\gf}^2+\myfrac {1}{\sin^2\gf}\abs{\nabla_{\gs'}v}^2\right)^{(p-2)/2}
=(\sin\gf)^{(\gb-1)(p-2)}\left(\gb^2w^2+\abs{\nabla_{\gs'}w}^2\right)^{(p-2)/2},
$$
thus
$$\BA{l}
\myfrac {1}{\sin^{N-2}\gf}\left[\sin^{N-2}\gf\left(\gb^2v^2+v_{\gf}^2+\myfrac {1}{\sin^2\gf}\abs{\nabla_{\gs'}v}^2\right)^{(p-2)/2}\!\!v_{\gf}\right]_{\gf}\\[6mm]
=\gb (\sin\gf)^{(\gb-1)(p-1)-1}
\left(\left(N-2+(\gb-1)(p-1)\right)-\left(N-1+(\gb-1)(p-1)\right)\sin^2\gf\right)\\[4mm]
{\ti}\left(\gb^2w^2+\abs{\nabla_{\gs'}w}^2\right)^{(p-2)/2}w,
\EA
$$
and
$$\BA {l}
\myfrac {1}{\sin^2\gf}div_{\gs'}\left[\left(\gb^2v^2+v_{\gf}^2+\myfrac {1}{\sin^2\gf}\abs{\nabla_{\gs'}v}^2\right)^{(p-2)/2}\!\!\nabla_{\gs'}v\right]\\[6mm]
\phantom {---}
=(\sin\gf)^{(\gb-1)(p-1)-1}div_{\gs'}\left[\left(\gb^2w^2+\abs{\nabla_{\gs'}w}^2\right)^{(p-2)/2}\nabla_{\gs'}w\right]
\EA$$
Finally $w$ satisfies
\begin {equation}\label {psrad-N-1}
-div_{\gs'}\left[\left(\gb^2w^2+\abs{\nabla_{\gs'}w}^2\right)^{(p-2)/2}\nabla_{\gs'}w\right]
=\gl_{N-1,\gb}\left(\gb^2w^2+\abs{\nabla_{\gs'}w}^2\right)^{(p-2)/2}w
\end {equation}
on $S^{N-2}$, which is the desired induction. \qeda
\\

In order to be more precise, we can completely represent the preceding solutions by introducing the generalized Euler angles in $\BBR^N=\{x=(x_{1},...,x_{N})\}$
\begin{equation}\label{Euler}\left\{\BA {l}
x_{1}=r\sin\gth_{N-1}\sin\gth_{N-2}...\sin\gth_{2}\sin\gth_{1}\\
x_{2}=r\sin\gth_{N-1}\sin\gth_{N-2}...\sin\gth_{2}\cos\gth_{1}\\[-1mm]
.\\[-1mm]
.\\[-1mm]
.\\[-1mm]
x_{N-1}=r\sin\gth_{N-1}\cos\gth_{N-2}\\
x_{N}=r\cos\gth_{N-1}
\EA\right.\end{equation}
where $\gth_{1}\in [0,2\gp]$ and $\gth_{k}\in [0,\gp]$, for $k=2,...,N-1$. Notice that $\gth_{N-1}$ is the variable $\gf$ in the representation (\ref {repres0}). The above theorem combined with the induction process yields  to the following.
\bth {repres2} Let $N>3$ and $p>1$. For any positive integer $k$ there exists $p$-harmonic functions
in $\BBR^N$ under the form
\begin{equation}\label{psrad-N1}
u(x)=(r\sin\gth_{N-1}\sin\gth_{N-2}...\sin\gth_{2})^{\gb_{k}}\gw_{k}(\gth_1)
\end{equation}
where $(\gb_{k},\gw_{k})$ are obtained in \rth {Krth}. Furthermore, if $p=N$ there exists a singular $N$-harmonic function under the form
\begin{equation}\label{psrad-N2}
u(x)=r^{-\gb_{k}}(\sin\gth_{N-1}\sin\gth_{N-2}...\sin\gth_{2})^{\gb_{k}}\gw_{k}(\gth_1).
\end{equation}
\es
\small{
Laboratoire de Math\'ematiques et Physique Th\'eorique\\
CNRS UMR 6083\\
Facult\'e des Sciences\\
Universit\'e Fran\c{c}ois Rabelais\\
F37200 Tours  France\\

\noindent borghol@univ-tours.fr\\
veronl@lmpt.univ-tours.fr}

\begin{thebibliography}{99}
\bibitem {Bo} Borghol R.,\textit{ Singularit\'es au bord de solutions d'\'equations quasilin\'eaires}, Th\`ese de Doctorat, Univ. Tours (in preparation).

\bibitem {BBV} Bidaut-V\'eron M. F., Borghol R. \& V\'eron L.,\textit{ Boundary Harnack inequalities and a priori estimates of singular solutions of quasilinear equations}, Calc. Var. and P. D. E., to appear.

\bibitem {FV} Friedman A., \& V\'eron L.,\textit{ Singular solutions of some quasilinear elliptic equations}, Arch. Rat. Mech. Anal. {\bf 96}, 359-387 (1986).


\bibitem{KV} Kichenassamy S. \& V\'eron L.,\textit{ Singular solutions of
the $p$-Laplace equation}, Math. Ann.  {\bf 275}, 599-615 (1986).

\bibitem{Kr} Krol I. N.,\textit{ The behavior of the solutions of a certain quasilinear equation near zero cusps of the boundary}, Proc. Steklov Inst. Math.  {\bf125}, 130-136 (1973).

\bibitem {Lie} Libermann G,\textit{ Boundary regularity for solutions of degenerate elliptic equations}, Nonlinear Anal.  {\bf12}, 1203-1219 (1988).

\bibitem{MW} Manfredi J. \& Weitsman A.,\textit{ On the Fatou Theorem for p-Harmonic Functions}, Comm. P. D. E.  {\bf13}, 651-668 (1988).

\bibitem {Re}Re\v setnjak, Ju. \textit{ Spatial mappings with bounded distortion} (Russian), 
Sibirsk. Mat. \v Z. {\bf 8}, 629-658 (1967) .

\bibitem{Se1} Serrin J.,\textit{ Local behaviour of solutions of
quasilinear equations}, Acta Math.  {\bf111}, 247-302 (1964).

\bibitem{SZ} Serrin J. \& Zou H.,\textit{ Cauchy-Liouville and universal boundedness theorems for
quasilinear elliptic equations and inequalities}, Acta Math.  {\bf189}, 79-142 (2002).

\bibitem{To} Tolksdorff P.,\textit{ On the Dirichlet problem for
quasilinear equations in domains with conical boundary points}, Comm.
Part. Diff. Equ.  {\bf 8}, 773-817 (1983).

\bibitem{To1} Tolksdorff P.,\textit{ Regularity for a more general class of quasilinear elliptic equations}, J. Diff. Equ.  {\bf 51}, 126-140 (1984).

\bibitem{Tr} Trudinger N.,\textit{ On Harnack type inequalities and
their applications to quasilinear elliptic equations}, Comm.
Pure Appl. Math.  {\bf 20}, 721-747 (1967).

\bibitem{Ve2} V\'eron L.,\textit{ Some existence and uniqueness
results for solution of some quasilinear elliptic equations on
compact Riemannian manifolds}, Colloquia Mathematica Societatis
J\'anos Bolyai  {\bf 62}, 317-352 (1991).

\bibitem{Ve3} V\'eron L.,\textit{ Singularities
of solutions of second order quasilinear elliptic equations}, Pitman
Research Notes in Math.  {\bf 353}, Addison-Wesley- Longman (1996).

\bibitem{Ve4} V\'eron L.,\textit{ Singularities of some quasilinear equations}, Nonlinear diffusion equations and their equilibrium states, II (Berkeley, CA, 1986), 333-365, Math. Sci. Res. Inst. Publ., {\bf 13}, Springer, New York (1988).

\bibitem {Vi}Vilenkin N. \textit{ Fonctions sp\'eciales et th\'eorie de la repr\'esentation des groupes},
Dunod, Paris (1969).

\end {thebibliography}
\end {document}